# Gaussian Fluid Queue with Autocorrelated Input

Kerry Fendick
78 Valley Avenue, Highlands, NJ, USA 07732
fendick@att.net

## *Abstract*

This paper develops a generalization of Brownian motion with stationary, autocorrelated increments as a tractable model for problems in business and finance. We show that any real continuous Gaussian Markov process with stationary increments and smooth covariance function is characterized by three parameters quantifying drift, volatility, and autocorrelations. We model a queue as a functional of a process defined by those characteristics and derive its transient distribution conditional on its history.

# 1 Introduction

Stochastic processes with continuous sample paths and stationary increments are natural models for time-varying phenomena ranging from security prices to inventory levels. Analyzing such phenomena is often simplified through the additional assumption of independent increments. Consequently, Brownian motion, a stochastic process with continuous sample paths and stationary, independent increments, has been widely applied in business and finance. A convenient characteristic of Brownian motion is its dependence on only a single parameter, which can be interpreted as a measure of variability or volatility. In applications, a time-varying quantity that changes frequently is commonly modeled as a functional of a process that is the sum of a scaled Brownian motion and a deterministic drift.

Independent increments are at once the most useful and most restrictive feature of Brownian motion. Functionals of Brownian motion can emerge as limits of more complicated systems, in which case the independence of increments itself emerges from more intuitive assumptions, e.g., heavy-traffic conditions for queues. Then, a model based on Brownian motion is a useful tool for developing approximations and controls for the more complicated system. But, more generally, the assumption of independent increments runs counter to intuition that what has just happened may help predict what will comes next. It therefore calls for strong justification for use in any particular application.

In this paper, we develop a model with continuous sample paths and the tractability of Brownian motion, but with the flexibility to quantify the degree of dependence between increments. Motivation for this work includes the modeling of physical queues and inventories.

## 1.1 Discussion of Results

Brownian motion is a Gaussian Markov process. Its Gaussian and Markov properties can be derived from its definition as a stochastic process with zero mean, continuous sample paths, and stationary, independent increments. In generalizing Brownian motion, we retain its Gaussian and Markov properties, while replacing the assumption that increments are stationary and independent with the weaker assumption that they are just stationary. We call this generalization a Gaussian Markov Stationary Increment (GMSI) process and show that it is characterized in great generality by only two parameters. A GMSI process with drift is therefore characterized by three parameters. Depending on these parameters, the increments of the process can be positively correlated, negatively correlated, or uncorrelated. A GMSI process with non-zero autocorrelations exhibits behavior called volatility clustering frequently observed for asset prices.

A prominent application of Brownian motion is the analysis of queues or inventories. The evolving state of a physical inventory or financial account can be modeled as a non-negative process that is a functional of the cumulative supply and demand for the queue's contents. The difference between the cumulative supply and demand is called the queue's net input. The non-negative queue-length process can be expressed as the sum of the net input process and a cumulative increasing process representing the queue's lost potential output or lost demand. We show that the class of GMSI processes is closed under particular operations of superposition and conditioning. We apply those properties to model a queue for which the net input process is a GMSI process that itself comprises a superposition of independent GMSI components. We then derive the conditional distribution of the queue length given its observed history and the history of each of those components.

The ability to model a queue's net input as a superposition is useful in applications, because one often starts with descriptions or measurements of independent sources of supply and demand for the queue's content. When a history is known for each source individually, conditioning on it provides more information than conditioning on the aggregated state of their superposition. As a prototypical example, a computer manufacturer may desire to forecast the future content of a shared inventory of a particular type of memory chip used in several different products lines (each corresponding to a particular desk-top or lap-top computer). Given a finite history for the net input process associated with each product line and for the content of the shared inventory, Theorem 3 of this paper describes the conditional distribution of the inventory content at a given future time when each net input process has been modeled as a GMSI process. The conditional distribution from Theorem 3 can be used to calculate the conditional expectation function or conditional quantile function for the inventory content, which in turn can be used for the desired forecast.

## 1.2 Related Work

The work here applies Doob's [1] characterization of Gaussian Markov processes, as formalized by Hida [2]. Results in Hida [2] also extended those of Doob to further characterize Gaussian Markov processes that are stationary. A GMSI process, however, is not itself stationary. Our results extend those of Doob to further characterize Gaussian Markov processes with stationary increments.

Reflected Brownian Motion (RBM) emerges as a heavy traffic limit for a wide class of queuing processes; for background, see Whitt [3] and references sited there. Iglehart & Whitt [4], Fendick, Saksena, & Whitt [5], and Fendick, Saksena, & Whitt [6] showed how the simple covariance structure of Brownian motion can approximate those of much more complex net input processes to queues in the heavy-traffic limit. Extremely complex covariance structures, including those from superposing an arbitrary number of other complex processes, can be captured through a handful of parameters in the heavy-traffic limit. Like RBM, the GMSI model for a queue in this paper is a diffusion process with a regulating barrier. It aims to characterize even more complex dependences through one additional degree of freedom than allowed by standard Brownian motion. As with the aforementioned Brownian models, our queuing results cover the case in which a net input process comprises a superposition of other independent, but potentially complex processes.

The additional degree of freedom enables the GMSI model to quantify autocorrelations of the net input process, a characteristic not captured by models based on standard Brownian motion. For a GMSI process, autocorrelations result in a conditionally-induced change of drift affecting the future increments of the process. Previously, Fendick [7] studied examples in which a net input process based on standard Brownian motion exhibited a conditionally-induced changes of drift, but explaining past behavior rather than predicting future behavior.

Given the foundation in heavy-traffic limit theorems for queues, authors including Newell [8] and Harrison [9] have modeled complex queuing systems directly as *stochastic flow systems* – systems with net input processes modeled as multi-dimensional Brownian motion with drift and with regulating barriers constraining behavior at boundaries. Adopting this philosophy, Fendick & Whitt [10] showed how the parameters for such a Brownian approximation of a queue with finite buffers can be directly estimated using historic observations of the queue's behavior. We follow a similar philosophy here but using GMSI processes instead of standard Brownian motion.

The transient distribution of RBM was derived by Newell [8] through analytical methods and by Harrison [9] through probabilistic arguments. We obtain the transient distribution of RBM as a limit of results here for a GMSI model and follow elements of Harrison's proof in generalizing his results. A stationary distribution for a Reflected Brownian Bridge (RBB) was derived by Hajek [11] as a model for a queue with periodic arrivals. We derive here the transient distribution for RBB as another special case of the GMSI model and then obtain Hajek's stationary distribution as another limit. These results provide examples of how a GMSI process generalizes both a Brownian motion and a Brownian bridge. The transient distribution for RBB derived here appears to be new.

The flexibility of GMSI models to represent both positively and negatively correlated increments is an essential characteristic of a generic approach to the modeling of queues. An example of a queue for which the net input process has negatively correlated increments is the model with periodic arrivals studied in Hajek [11]. An example of a queue for which the net input process has positively correlated increments is the model with bursty arrivals studied in Fendick, Saksena, and Whitt [6] and interpreted in (48)-(51) on page 181 of Fendick and Whitt [12].

The queuing model studied here is an example of a Gaussian fluid queue, as defined and studied by other authors under variations on the assumptions made here. Debicki and Roski [13], Debicki, Es-Saghouani, and Mandjes [14], and Debicki, Koskinski, and Mandejes [15] derive asymptotic results for queues assuming that the net input process is Gaussian with stationary increments, but not assuming that it is Markov. Debicki and Roski [16] derive asymptotic results for queues assuming that the net input process is Gaussian Markov and stationary.

In a recent publication, Debicki and Mandjes [17] describe open problems in the theory of Gaussian fluid queues including speed of convergence to stationarity and correlation structure for the queueing process. They note that that "for a general Gaussian input process …. there are no explicit expressions available for the (stationary) distribution of (the queue length) … let alone the transient". For a GMSI net input process, however, Theorem 3 of the current paper provides an explicit expression for the transient queue length distribution as a special case. It therefore can be used to explore the approach to stationarity. Since Theorem 3 also describes the transient queue-length distribution at a given time conditional on the of queue length at a given earlier time, the joint distribution of the queue lengths at those two times – and hence the correlation structure of the queueing process - can in principle be obtained as a corollaries from the definition of a conditional distribution.

In one sequel to this paper, Fendick [18] describes methods for estimating parameters for a multidimensional GMSI model. The methods there cover the use of samples from a given process and from an arbitrary number of other processes on which the given process is statistically dependent. They also cover parameter estimation for the queuing applications of this paper. In particular, they cover estimation of the conditional expectation functions for the length of the queue, its net input, and lost potential output.

Fendick [18] shows how a GMSI process can be used to model the local behavior of a Gaussian process for which increments are stationary over all time scales, but Markov over only limited time scales. That modeling approach is generally useful to optimize the fit of model parameters for the time scales that matter most in a particular application. It is essential for applying GMSI processes when autocorrelations are negative, because the time domain of a GMSI process with negative autocorrelations is necessary a finite interval, as the results of Section 2 of the current paper imply.

In a second sequel to this paper, Fendick [19] studies the problem of rational option pricing for a model in which the log returns of the underlying security evolves according to a Gaussian Markov process with ergodic properties but non-stationary increments. It shows how GMSI processes characterize such log returns over large time scales.

## 2 Gaussian Markov Processes with Stationary Increments

Let $X(t)$ denote a real, zero-mean stochastic process on $0 \leq t < \Delta \leq \infty$ with $X(0) = 0$, and assume that

$$w(s,t) \equiv E\left[\left(X(t+s) - X(s)\right)^2\right] \tag{2.1}$$

is finite on $0 \leq s \leq s + t < \Delta$. When $w(s,t)$ does not vary with $s$, the increments of $X$ are said to be wide-sense stationary; see for example pages 99-101 of Doob [20]. When, in addition, $X$ is a Gaussian process, the increments are also strictly stationary, since the finite-dimensional distribution of a zero-mean Gaussian process is determined by its covariances.

In the statements that follow, we will let $B(t)$ for $0 \leq t < \infty$ denote a standard Brownian motion, which is the unique continuous, zero-mean Gaussian process with $B(0) = 0$ and covariance function

$$E[B(s)B(t)] = s \text{ for any } 0 \leq s \leq t. \tag{2.2}$$

We will use "smooth" to describe a function that is twice continuously differentiable over its entire domain in each of its arguments.

*Theorem 1*: *If* $\{X(t): 0 \leq t < \Delta\}$ *is a real, continuous, zero-mean Gaussian process with* $X(0) = 0$ *and smooth covariance function* $E[X(s)X(t)] = r(s,t)$ *on* $0 \leq s \leq t < \Delta \leq \infty$, *then X has stationary increments and the Markov property if and only if*

$$r(s,t) = s(\theta - \tau t) \tag{2.3}$$

*for some* $\theta > 0$ *and* $\tau \leq \theta/\Delta$.

Proof of Theorem 1:

We first assume that $X$ is a real, continuous, zero-mean Gaussian Markov process on $[0, \Delta)$ with stationary increments and smooth covariance function. Since $X$ is continuous, Gaussian, and Markov, Theorem II.1 of Hida [2] implies that there exist real valued functions $f(\cdot)$ and $g(\cdot)$ with $f(s) \neq 0$ for $0 \leq s < \Delta$ such that

$$X(t) = f(t) \int_0^t g(s)\, dB(s) \text{ for } 0 \leq t < \Delta. \tag{2.4}$$

It follows that $X$ has covariance function

$$r(s,t) = G(s)f(s)f(t) \text{ for } 0 \leq s \leq t < \Delta \tag{2.5}$$

where

$$G(s) \equiv \int_0^s g^2(u)\, du. \tag{2.6}$$

Since $r(\cdot,\cdot)$ is smooth, $f(\cdot)$ and $G(\cdot)$ are too, as implied by (2.5). By (2.1) and (2.5)

$$w(s,t) = G(t+s)f^2(t+s) - 2G(s)f(s)f(t+s) + G(s)f^2(s) \tag{2.7}$$

for $0 \leq s \leq s + t < \Delta$, so that $w(s,t)$ is also smooth. If $X(\cdot)$ has stationary increments, then we must also have $\partial/_{\partial s}\, w(s,t) = 0$ on $0 \leq s \leq s + t < \Delta$. In particular,

$$0 = \frac{\partial}{\partial s} w(0,t) = \left(f^2(t)G(t)\right)' + f(0)g^2(0)\left(f(0) - 2f(t)\right) \tag{2.8}$$

for $0 \leq t < \varDelta$, where the second equality of (2.8) follows from (2.5)-(2.7).

Therefore,

$$G(s) = \frac{f(0)g^2(0)}{f^2(s)}\left(2\int_0^s f(u)\,du - f(0)s\right) \qquad (2.9)$$

for $0 \leq s < \varDelta$. For a non-trivial solution, we will assume that $g(0) \neq 0$.

Likewise, we must have $\partial^2/\partial s^2\, w(s,t) = 0$ on $0 \leq s \leq s + t < \varDelta$. Using (2.9), we find that the terms of $\partial^2/\partial s^2\, w(s,t)$ involving $f''$ vanish as $s$ approaches the origin, and we thus obtain

$$0 = \frac{\partial^2}{\partial s^2} w(0,t) = 2f(0)g^2(0)\big(f'(0) - f'(t)\big) \qquad (2.10)$$

for $0 \leq t < \varDelta$. Since $f(0) \neq 0$ and $g(0) \neq 0$, we conclude that

$$f(t) = f(0) + f'(0)t \; for \; 0 \leq t < \Delta. \qquad (2.11)$$

From (2.5), (2.9), and (2.11), we arrive at (2.3) for some constants θ and τ. The condition that $E[X^2(t)] \geq 0$ for $0 \leq t < \Delta$ implies that $\theta > 0$ and $\tau \leq \theta/\Delta$ in (2.3).

Conversely, if $r(s,t) = s(\theta - \tau t)$, then $r(\cdot,\cdot)$ is smooth. If $X$ is a real, continuous, zero-mean Gaussian process satisfying (2.3), then $w(s,t)$ in (2.1) does not vary with $s$, so that $X$ has stationary increments; and we verify that $r(\cdot,\cdot)$ satisfies the conditions of Theorem 8.1 on p. 233 of Doob [20] under which $X$ is a Markov process in the wide sense. Since $X$ is real and Gaussian, the wide-sense Markov property of $X$ implies that $X$ is also a Markov process in the strict sense; see the remark at the bottom of p. 233 of Doob [20]. ■

We will say that $X$ is a $(\theta,\tau)$ $GMSI$ process on $[0,\varDelta)$ if $X$ is a real, continuous, zero-mean Gaussian Markov process with stationary increments and smooth covariance function characterized by the given parameters. When $X$ is a $(\theta,\tau)$ $GMSI$ process on $[0,\varDelta)$ and $\rho$ is an arbitrary constant, the process $\rho t + X(t)$ for $0 \leq t < \Delta$ is also Gaussian Markov with stationary increments and is characterized by the triple $(\rho,\theta,\tau)$.

*Corollary 1: When* $\{B(t): t \geq 0\}$ *denotes standard Brownian motion, the following statements are equivalent:*

i. $X$ is a $(\theta,\tau)$ $GMSI$ process on $[0,\varDelta)$.

ii. $X(t) = (\theta - \tau t)B\left(\frac{t}{\theta - \tau t}\right)$ for $0 \leq t < \varDelta$.

iii. *X is a diffusion process such that* $X(0) = 0$,

$lim_{h\to 0}\, h^{-1}\, E[X(t+h) - X(t) | X(t) = x] = -\tau x/(\theta - \tau t)$, *and*

$lim_{h\to 0}\, h^{-1}\, E\left[\left(X(t+h) - X(t)\right)^2 | X(t) = x\right] = \theta.$

Proof of Corollary 1: The equivalence of (i) and (ii) can be verified through (2.2) and (2.3). By (2.3)-(2.6), which hold under (i),

$$dX(t) = \frac{-\tau}{\theta - \tau t} X(t)dt + \sqrt{\theta} dB(t)\ for\ 0 \le t < \Delta, \tag{2.12}$$

and the diffusion coefficients in (iii) can be read from (2.12) as described in Chapter 15, Section 14 of Karlin and Taylor [21]. ■

For random variables $X$ and $Y$, let $Cor(X,Y) \equiv Cov(X,Y)/\left(\sqrt{Var\ X}\sqrt{Var\ Y}\right)$. As is well known, $-1 \le Cor(X,Y) \le 1$ whenever it is finite.

*Corollary 2: Let X denote a* $(\theta, \tau)$ *GMSI process on* $[0, \Delta)$. *Assume that* $t \ge 0$, $u \ge s \ge 0$, *and* $t + u + s < \Delta$, *and let* $R_s(t) \equiv X(t+s) - X(t)$. *Then*

$$Cor\left(R_s(t), R_s(t+u)\right) = \frac{-s\tau}{\theta - s\tau}\ and\ Cor\left(R_s^2(t), R_s^2(t+u)\right) = \frac{s^2\tau^2}{(\Theta - sT)^2}, \tag{2.13}$$

*as the result of which*

$$\lim_{\tau\to -\infty} Cor\left(R_s(t), R_s(t+u)\right) = 1,\quad \lim_{\substack{\tau\to\theta/\Delta \\ s\to\Delta/2}} Cor\left(R_s(t), R_s(t+u)\right) = -1,$$

and

$$\left(Cor\left(R_s(t), R_s(t+u)\right)\right)^2 \le Cor\left(R_s^2(t), R_s^2(t+u)\right) \le \left|Cor\left(R_s(t), R_s(t+u)\right)\right|. \tag{2.14}$$

Proof of Corollary 2: The first equality of (2.13) follows immediately from (2.3). The second equality of (2.13) also follows from (2.3) since Gaussian moments of any order reduce to covariances; c.f. Isselis [22]. ■

The quantities $Cor\left(R_s(t), R_s(t+u)\right)$ are commonly called autocorrelation coefficients for time scale $s$. As the first statement of (2.13) implies, a $(\theta, \tau)$ $GMSI$ process has negative autocorrelation coefficients when $\tau$ is positive, and positive autocorrelation coefficients when $\tau$ is negative.

The second statement of (2.13) implies that periods of high volatility (variability) will tend to follow periods of high volatility and periods of low volatility will tend to follow periods of low volatility. Such behavior is called volatility clustering in the financial literature. By (2.14), volatility clustering disappears for a GMSI process as autocorrelation coefficients approaches zero.

Corollary 2 shows that a given GMSI process always has autocorrelations of the same sign. Theorem 1 implies that the time domain of a GMSI process $X$ must be bounded if $X$ has negative autocorrelations but may be unbounded if $X$ has positive autocorrelations.

## 3 An Example

Let $B(\cdot)$ denote a standard Brownian motion, as defined in Section 2. Then, $B(\cdot)$ satisfies the conditions of Theorem 1 with $\theta = 1$ and $\tau = 0$, so that $B$ is a $(1,0)$ $GMSI$ process on $[0,\infty)$. Next, let $B_\Delta^0(t)$ denote a Brownian bridge on $0 \leq t < \Delta$, defined here as a continuous zero-mean Gaussian process with $B_\Delta^0(0) = 0$ and covariance function

$$E\left[B_\Delta^0(s)B_\Delta^0(t)\right] = \frac{s}{\Delta}\left(1 - \frac{t}{\Delta}\right) \text{ on } 0 \leq s \leq t \leq \Delta. \tag{3.1}$$

It follows that $B_\Delta^0$ is a $(\Delta^{-1}, \Delta^{-2})$ $GMSI$ process on $[0, \Delta)$.

By (3.1), $E\left[\left(B_\Delta^0(\Delta)\right)^2\right] = 0$, from which it follows that $B_\Delta^0(\Delta) - B_\Delta^0(s) = -B_\Delta^0(s)$ with probability one. In other words, $B_\Delta^0(\cdot)'s$ increment on the interval $[s, \Delta]$ is always equal in absolute value but opposite in sign to its increment on $[0, s]$. This is one extreme in dependence structure.

The Brownian motion $B(\cdot)$ exhibits another extreme since a Brownian motion has independent increments; see for example Harrison [9]. The independence of increments is reflected by the covariance function (2.2), which implies that $E[B(s)\big(B(t) - B(s)\big)] = 0$ for any $0 \leq s \leq t$.

Let

$$X(t) \equiv (1 - \chi)\sqrt{\alpha}B_\Delta^0(t) + \chi\sqrt{\beta}B(t) \text{ for } 0 \leq t < \Delta, \tag{3.2}$$

where $0 \leq \chi \leq 1$, $\alpha > 0$, and $\beta > 0$ are constants, and where $B(\cdot)$ and $B_\Delta^0(\cdot)$ are assumed to be independent of one another.

The following theorem follows from Theorem 1 as is easily verified by comparing the relevant covariance functions.

*Theorem 2*: *If X is given by (3.2), then X is a* $(\theta, \tau)$ *GMSI process on* $[0, \Delta)$ *where*

$$\theta \equiv \frac{\alpha(1-\chi)^2 + \Delta\beta\chi^2}{\Delta} > 0 \text{ and } 0 \leq \tau \equiv \frac{\alpha(1-\chi)^2}{\Delta^2} \leq \frac{\theta}{\Delta}.$$

*Conversely, if* $X$ *is a* $(\theta,\tau)$ *GMSI process on* $[0,\Delta)$ *where* $0 \le \tau \le \theta/\Delta$, *then* $X(t) = \sqrt{\tau\Delta^2}B_\Delta^0(t) + \sqrt{\theta - \tau\Delta}B(t)$ *is a representation for* $0 \le t < \Delta$.

By Corollary 2 and Theorem 2, (3.2) is a canonical representation of a $(\theta,\tau)$ $GMSI$ process for which autocorrelations are negative.

When $X$ is a $(\theta,\tau)$ $GMSI$ process on $[0,\Delta)$ where $\Delta < \infty$ and $\tau = \theta/\Delta > 0$, it follows from (2.3) that $\lim_{s\to\Delta} EX^2(s) = 0$. We will then say that $X$ is a scaled Brownian bridge.

## 4 Closure Properties of GMSI Processes

The following propositions describe closure properties of the class of GMSI processes under particular operations of superposition and conditioning. These propositions are used for the queuing analysis in Section 5 and are useful in applications more generally.

For $i = 1,\dots,K$, let $\theta^{(i)} > 0$ and $\tau^{(i)} \neq 0$ denote constants, and assume that $X^{(i)}$ is a $\left(\theta^{(i)},\tau^{(i)}\right)$ $GMSI$ *process on* $[0,\Delta)$ where $0 < \Delta \le \min_{0\le i\le K}\theta^{(i)}/\tau^{(i)}$ if $\tau^{(i)} > 0$ for some $i \in 1,\dots,K$; and $0 < \Delta \le \infty$ otherwise. Assume further that the $X^{(i)}{}'s$ are independent of one another, and let

$$X(t) \equiv \sum_{i=1}^{K} k_i X^{(i)}(t) \text{ on } 0 \le t \le \Delta \tag{4.1}$$

*where the* $k_i{}'s$ *are constants.*

Our first proposition follows from Theorem 1 and the linear product form of the covariance function of a GMSI process.

*Proposition 1: If* $X$ *is defined as in* $(4.1)$, *then it is a* $(\theta,\tau)$ *GMSI process on* $[0,\Delta)$ *where*

$$\theta \equiv \sum_{i=1}^{K} k_i^2\theta^{(i)} > 0 \quad and \quad \tau \equiv \sum_{i=1}^{K} k_i^2\tau^{(i)} \le \frac{\theta}{\Delta}. \tag{4.2}$$

Proposition 1 shows that the Gaussian and Markov properties and independence of increments are preserved under superpositions of independent GMSI processes. Preservation of the Markov process under such superpositions is a trait of GMSI processes not shared with processes more generally.

Proofs of the following two propositions can be obtained by finding diffusion coefficients, as exemplified on page 269 of Karlin and Taylor [21], for the conditioned processes below, and comparing them to the diffusion coefficients for a GMSI process in (iii) of Corollary 1.

For $X(\cdot)$ defined in (4.1), $0 \le u < \Delta$, and $0 \le h < \Delta - u$, let

$$X_{u;\vec{x}}(h) = X(u+h) - X(u) + h\sum_{i=1}^{K} \frac{k_i x_i}{\theta^{(i)}/\tau^{(i)} - u} \textit{ conditional on } X^{(i)}(u) = x_i \textit{ for } i = 1, \dots, K\,. \tag{4.3}$$

In other words,

$$\{X_{u;\vec{x}}(h): 0 \le h < \Delta - u\} = \left\{X(u+h) - X(u) + h\sum_{i=1}^{K} \frac{k_i x_i}{\theta^{(i)}/\tau^{(i)} - u} : 0 \le h < \Delta - u\right\}$$

confined to sample paths for which $X^{(i)}(u) = x_i$ for $i = 1, \dots, K$.

Proposition 2. *If $X_{u;\vec{x}}$ is defined as in (4.3), then it is a $(\theta, \tau_u)$ GMSI process on $[0, \Delta - u)$ where $\theta$ is defined in (4.2) and*

$$\tau_u \equiv \sum_{i=1}^{K} k_i^2 \frac{\theta^{(i)}}{\theta^{(i)}/\tau^{(i)} - u} \le \frac{\theta}{\Delta - u}. \tag{4.4}$$

From Proposition 2 and the definition of $X_{u;\vec{x}}$ in (4.3), we see that conditioning on $X^{(i)}(u) = x_i$ for $i = 1, \dots, K$ induces $X(s)$ to have a constant drift of $\sum_{i=1}^{K} \frac{-k_i x_i}{\theta^{(i)}/\tau^{(i)} - u}$ on $u \le s < \Delta$. The special case of conditioning on $X(u) = x_1$ alone is obtained when $K = 1$ and $k_1 = 1$. For that case, the conditional mean of $X(s)$ on $u \le s < \Delta$ reverts towards zero if $X(\cdot)$ has negative autocorrelations and away from zero if $X(\cdot)$ has positive autocorrelations.

For $X$ still defined in (4.1), next let

$$Z(t) \equiv \rho t + X(t) \text{ for } 0 \le t < \Delta \tag{4.5}$$

where $\rho$ is arbitrary constant. Then, for $0 \le u < \Delta$ and $0 < h < w < \Delta - u$, let

$$Z_{u;\vec{x}}^{u+w;z}(h) \equiv Z(u+h) - Z(u) - (z/w)h \textit{ conditional on } X^{(i)}(u) = x_i \textit{ for } i = 1, \dots, K \textit{ and on } Z(u+w) - Z(u) = z\,. \tag{4.6}$$

*Proposition 3: If $Z_{u;\vec{x}}^{u+w;z}$ is defined as in (4.6), then it is a $(\theta, \theta/w)$ GMSI process on $[0, w)$ where $\theta$ is defined by (4.2).*

Similarly to Proposition 2, the induced drift of $Z(s)$ on $u \le s < u + w$ is constant. But, in this case, it depends neither on the values $x_i$ for $i = 1, \dots, K$ nor on the original drift $\rho$. As an explanation, note that the process $Z_{u;\vec{x}}^{u+w;z}(h)$ is a scaled Brownian bridge on $0 \le h < w$ and so approaches zero with probability one as $h$ approaches $w$. To account for the increment $Z(u+w) - Z(u) = z$, the induced drift of $Z$ over the interval of length $w$ must therefore equal $z/w$.

## 5 Queues Driven By GMSI Processes

Recalling that $Z$ in (4.5) depends on parameters $\rho, \{\theta^{(i)}\}_{i=1}^{K}$, and $\{\tau^{(i)}\}_{i=1}^{K}$ through $X$ as defined in (4.1), let

$$Q(t) = Q(0) + Z(t) + L(t) \tag{5.1}$$

where $Q(0) \geq 0$, and $L(\cdot)$ is a non-decreasing continuous function with the properties that $L(0) = 0$ and that $L(\cdot)$ increases only when $Q(\cdot) = 0$ such that $Q(t) \geq 0$ for all $0 \leq t \leq \Delta$. Then, $Q(t)$ has the interpretation of a non-negative queue length at time $t$ with an initial value of $Q(0)$. The function $L(\cdot)$ has the interpretation of the queue's lost potential output, and $Z(\cdot)$ of the queue's net input process equal to the difference between the queue's cumulative input (or supply) and its cumulative potential output (or demand), including any lost potential output.

As an example, if $Q(t)$ models an inventory level at time $t$, then $Z(t)$ represents the value of new inventory purchased or produced over the interval $[0, t]$ minus the value of demand for inventory over that same period. And $L(t)$ represents the value of additional inventory that would have been sold over the interval $[0, t]$ if the firm had maintained sufficient inventory levels to meet all demand over that interval.

For $Q$ defined in (5.1), $0 \leq u < \Delta,\ 0 < h < \Delta - u$, let

$$Q_{u;\vec{x},v}(h) \equiv Q(u+h) \ \textit{conditional on}\ Q(u) = v\ \textit{and on}\ X^{(i)}(u) = x_i\ \textit{for}\ i = 1, \ldots, K. \tag{5.2}$$

*Theorem 3: If $Q_{u;\vec{x},v}$ is defined as in (5.2), then*

$$\begin{aligned} F_{u;\vec{x},v}(h;q) &\equiv P\left(Q_{u;\vec{x},v}(h) \leq q\right) \\ &= \frac{1}{2}\left(1 - e^{\frac{-2q(\tau_u q - \theta\rho_{u;\vec{x}})}{\theta^2}} - erf\left(\frac{-q + v + \rho_{u;\vec{x}}h}{\sqrt{2h(\theta - \tau_u h)}}\right)\right. \\ &\quad \left. + e^{\frac{-2q(\tau_u q - \theta\rho_{u;\vec{x}})}{\theta^2}} erf\left(\frac{\theta(q+v) - (2\tau_u q - \theta\rho_{u;\vec{x}})h}{\theta\sqrt{2h(\theta - \tau_u h)}}\right)\right) \end{aligned} \tag{5.3}$$

*where* $erf(z) \equiv 2\pi^{-1/2}\int_0^z \exp(-t^2)\,dt,$ $\theta$ *is defined in (4.2),* $\tau_u$ *is defined in* (4.4), *and*

$$\rho_{u;\vec{x}} \equiv \rho - \sum_{i=1}^{K} \frac{k_i x_i}{\theta^{(i)}/\tau^{(i)} - u}. \tag{5.4}$$

*Proof of Theorem 3*:

We first prove Theorem 3 for $0 \le \tau_u \le \theta/(\Delta - u)$.

Because $X(\cdot)$ and $X^{(i)}(\cdot)$ for $i = 1, \dots, K$ are Markov, as follows from Theorem 1, $Q_{u;\vec{x},v}$ is independent of the states of $Q(t)$ and $X^{(i)}(t)$ for $i = 1, \dots, K$ on $0 \le t < u$. Let

$$V_u(h) \equiv Z(u+h) - Z(u) \tag{5.5}$$

and

$$V_{u;\vec{x},v}(h) \equiv V_u(h) \textit{ conditional on } X^{(i)}(u) = x_i \text{ for } i = 1, \dots, K, \textit{ and } Q(u) = v. \tag{5.6}$$

For $Q_{u;\vec{x},v}$ defined in (5.2) and $V_{u;\vec{x},v}$ in (5.6), we have

$$Q_{u;\vec{x},v}(h) = v + V_{u;\vec{x},v}(h) + L_{u;\vec{x},v}(h)$$

where $L_{u;\vec{x},v}(0) = 0$ and $L_{u;\vec{x},v}(\cdot)$ increases only when $Q_{u;\vec{x},v}(\cdot) = 0$ such that $Q_{u;\vec{x},v}(h) \ge 0$ for all $0 \le h < \Delta - u$. It follows from Chaper 2, Section2, Propositions (10) of Harrison [9] that

$$Q_{u;\vec{x},v}(h) = \begin{cases} v + V_{u;\vec{x},v}(h), & \inf\limits_{0 \le s \le h} V_{u;\vec{x},v}(s) > -v \\ v + V_{u;\vec{x},v}(h) - \inf\limits_{0 \le s \le h} V_{u;\vec{x},v}(s), & \textit{otherwise} \end{cases} \tag{5.7}$$

Rewriting (5.5) as

$$V_u(h) = \rho h - \left(\sum_{i=1}^{K} \frac{k_i x_i}{\theta^{(i)}/\tau^{(i)} - u}\right) h + X(u+h) - X(u) + h \sum_{i=1}^{K} \frac{k_i x_i}{\theta^{(i)}/\tau^{(i)} - u}$$

as follows from (4.1) and (4.5), we see from (5.6) that

$$V_{u;\vec{x},v}(h) = \rho_{u;\vec{x}} h + X_{u;\vec{x}}(h) \; for \; 0 \le h < \Delta - u \tag{5.8}$$

where $X_{u;\vec{x}}$ is defined in (4.3) and $\rho_{u;\vec{x}}$ in (5.4).

Let

$$B(t) = \frac{1 + t\tau_u}{\theta} X_{u;\vec{x}}\left(\frac{t\theta}{1 + t\tau_u}\right) \text{ for } t \ge 0. \tag{5.9}$$

which, as Corollary 1 implies, is standard Brownian motion. The argument of $X_{u;\vec{x}}(\cdot)$ in (5.9) is an increasing function $t$. Because $0 \le \tau_u \le \theta/(\Delta - u)$, the argument of $X_{u;\vec{x}}(\cdot)$ in (5.9) is constrained to the interval $[0, \Delta - u)$, and the term $\frac{1 + t\tau_u}{\theta}$ multiplying $X_{u;\vec{x}}(\cdot)$ in (5.9) is always positive. Then,

$$\begin{aligned} Y(t) &\equiv v + V_{u;\vec{x},v}\left(\frac{t\theta}{1 + t\tau_u}\right) \\ &= v + \rho_{u;\vec{x}} \frac{t\theta}{1 + t\tau_u} + X_{u;\vec{x}}\left(\frac{t\theta}{1 + t\tau_u}\right) \\ &= v + \rho_{u;\vec{x}} \frac{t\theta}{1 + t\tau_u} + \frac{\theta}{1 + t\tau_u} B(t) \end{aligned} \tag{5.10}$$

using (5.8) and (5.9).

As an intermediate step, we will find

$$G(x,y) \equiv P\left(Y(t) \le x,\ \inf_{0\le s\le t} Y(s) \le y\right). \tag{5.11}$$

Because $1 + s\tau_u > 0$ whenever $s \ge 0$, we see from (5.10) that

$$\begin{aligned}
\inf_{0\le s\le t} Y(s) \le y &\Leftrightarrow \inf_{0\le s\le t}\left(v + \rho_{u;\vec{x}}\frac{s\theta}{1+s\tau_u} + \frac{\theta}{1+s\tau_u}B(t)\right) \le y \\
&\Leftrightarrow \inf_{0\le s\le t}\left(\frac{\rho_{u;\vec{x}}\theta s + \theta B(s) + (v-y)(1+s\tau_u)}{1+s\tau_u}\right) \le 0 \\
&\Leftrightarrow \inf_{0\le s\le t}\left(\rho_{u;\vec{x}}\theta s + \theta B(s) + (v-y)(1+s\tau_u)\right) \le 0 \\
&\Leftrightarrow \inf_{0\le s\le t}\left(\mu s + \theta B(s)\right) \le y - v
\end{aligned} \tag{5.12}$$

where $\mu \equiv \rho_{u;\vec{x}}\theta + (v-y)\tau_u$.

By (5.10)-(5.12),

$$G(x,y) = P\left(\mu t + \theta B(t) \le (x-y)(1+t\tau_u) + y - v,\ \inf_{0\le s\le t}\left(\mu s + \theta B(s)\right) \le y - v\right) \tag{5.13}$$

so that we can find $G(x,y)$ based solely on properties of $\{\mu s + \theta B(s): s \ge 0\}$, a scaled standard Brownian motion with constant drift. Applying Proposition 2 from Chapter 1, Section 8 of Harrison [9] and obtaining limits of integration from (5.13), we then see that

$$\begin{aligned}
G(x,y) &= \int_{-\infty}^{y-v}\int_{b}^{(x-y)(1+t\tau_u)+y-v} s(a,b,\mu)\,da\,db \\
&= \int_{-\infty}^{y-v}\int_{-\infty}^{a} s(a,b,\mu)\,db\,da + \int_{y-v}^{(x-y)(1+t\tau_u)+y-v}\int_{-\infty}^{y-v} s(a,b,\mu)\,db\,da
\end{aligned}$$

where $s(a,b,\mu) \equiv exp\left(\frac{\mu a}{\theta^2} - \frac{\mu^2 t}{2\theta^2}\right)\frac{d}{db}\phi\left(\frac{2b-a}{\theta t^{1/2}}\right)\frac{1}{\theta t^{1/2}}$ for $\phi(w) \equiv \exp\left(-\frac{w^2}{2}\right)/\sqrt{2\pi}$.

Then,

$$\begin{aligned}
g(x,y) &\equiv \frac{d}{dy}\frac{d}{dx}G(x,y) \\
&= (1+t\tau_u)\frac{d}{dy}\int_{-\infty}^{y-v} s\left((x-y)(1+t\tau_u) + y - v, b, \rho_{u;x}\theta + (v-y)\tau_u\right)db \\
&= \frac{\sqrt{2}(v+x-2y)(1+t\tau_u)^2 e^{-\left(\frac{t\rho_{u;\vec{x}}^2}{2} + \frac{(1+t\tau_u)(v-x)\rho_{u;\vec{x}}}{\theta} + \frac{(1+t\tau_u)\left((1+\tau_u)(x^2+v^2)+4y(y-x-v)+2vx(1-t\tau_u)\right)}{2t\theta^2}\right)}}{\sqrt{\pi}t^{3/2}\theta^3}
\end{aligned}$$

Using the above definition of $g(\cdot,\cdot)$ along with (5.7), (5.10), and (5.11), we find that

$$\begin{aligned}
F_{u;\vec{x},v}\left(\frac{t\theta}{1+t\tau_u};q\right) &= P\left(Y(t) - \inf_{0\le s\le t} Y(s) \le q,\ \inf_{0\le s\le t} Y(s) \le 0\right) \\
&\quad + P\left(Y(t) \le q,\ \inf_{0\le s\le t} Y(s) > 0\right) \\
&= \int_{-\infty}^{0}\int_{y}^{q+y} g(x,y)\,dx\,dy + \int_{0}^{v}\int_{y}^{q} g(x,y)\,dx\,dy.
\end{aligned}$$

so that

$$\frac{d}{dq}F_{u;\vec{x},v}\left(\frac{t\theta}{1+t\tau_u};q\right)=\int_{-\infty}^{0}g(q+y,y)dy+\int_{0}^{v}g(q,y)dy\,.$$

This last expression is a probability density function. The integrals on the right-hand side can be solved using above explicit expression for $g(\cdot,\cdot)$ by completing the squares in the exponent, and it is then trivial to verify agreement of the result with the corresponding density obtained from (5.3). A related example of completing the squares in the exponent is given on page 13 of Harrison [9]. Considering now the case in which $\tau_u < 0$, we replace the condition in (5.9) that $t \geq 0$ with the condition that $0 \leq t < (\Delta - u)/(\theta - (\Delta - u)\tau_u)$. Then, the argument of $X_{u;\vec{x}}(\cdot)$ in (5.9) is still constrained to the interval $[0, \Delta - u)$, and the term $\frac{1+t\tau_u}{\theta}$ multiplying $X_{u;\vec{x}}(\cdot)$ in (5.9) is still always positive. With this modification, the remainder of the proof holds with no additional changes. ■

The (unconditional) distribution of $Q(t)$ with initial conditions

$$X^{(i)}(0)=0 \text{ for } i=1,\ldots,K, \text{ and } Q(0)=v$$

is given by $F_{0;\vec{0},v}(t;q)$. The formula for $F_{0;\vec{0},v}(t;q)$ is given by the right-hand side of equation (5.3) with $\rho$ replacing $\rho_{u;\vec{x}}$, with $\tau$ from (4.2) replacing $\tau_u$, and with $t$ replacing $h$. When $X$ in (4.1) is a scaled Brownian bridge, corresponding to the case in which $\tau > 0$ and $\Delta = \theta/\tau$, the process $Q(\cdot)$ in (5.1) may be called a Regulated Brownian Bridge (RBB). We obtain a stationary distribution for RBB by taking the limit of $F_{0;\vec{0},v}(t;\cdot)$ as $t$ approaches $\theta/\tau$:

*Corollary 3 (stationary distribution for RBB): For* $F_{u;\vec{x},v}(\cdot;\cdot)$ *defined as in (5.3),* $\rho \equiv \rho_{0;\vec{0}}$ , *and* $\tau \equiv \tau_0 > 0$,

$$\vec{F}(q) \equiv \lim_{t\to\theta/\tau} F_{0;\vec{0},v}(t;q) = 1-e^{-2q(\tau q-\theta\rho)/\theta^2}. \qquad (5.14)$$

Previously, Hajek [11] derived the stationary distribution for RBB as a model for a queue with periodic arrivals, uniform phases, a constant service rate, and initial queue state of zero. For Hajek's model, $v = 0$, and the net input process corresponding to $Z(t)$ in (5.1) is $-\frac{\mathcal{M}-\mathcal{K}}{\mathcal{M}}t+\frac{\sqrt{\mathcal{K}}}{\mathcal{M}}B_1^0(t)$ , where $B_1^0(\cdot)$ is a zero-drift Brownian bridge with covariance function $E[B_1^0(s)B_1^0(t)] = s(1-t)$ on $0 \leq s \leq t \leq 1$. (There, $\mathcal{K}$ denotes the number of independent periodic sources, each generating $1/\mathcal{M}$ units of work once every time unit. ) To obtain the unconditional stationary distribution for this model, we use (5.14) with $\rho = -(\mathcal{M}-\mathcal{K})/\mathcal{M}$ and $\theta = \tau = \mathcal{K}/\mathcal{M}^2$ to obtain

$$\vec{F}(q) = 1 - exp(-2q\mathcal{M}(\mathcal{M} + \mathcal{M}q - \mathcal{K})/\mathcal{K}),$$

in agreement with (3.7) of Hajek [11].

For Regulated Brownian Motion (RBM), the net input process corresponding to $Z(\cdot)$ in (5.1) is $\rho t + \sqrt{\theta}B(t)$ where $B(\cdot)$ is a zero-drift Brownian motion with the covariance function given in (2.2). We therefore see that RBM corresponds to the limiting distribution of $Q(\cdot)$ in (5.1) as $\tau$ in (4.2) approaches zero.

*Corollary 4 (transient distribution for RBM): For* $F_{u;\vec{x},v}(\cdot;\cdot)$ *defined in (5.3),* $\rho \equiv \rho_{0;\vec{0}}$ , and $\tau \equiv \tau_0$,

$$\lim_{\tau\to 0} F_{0;\vec{0},v}(t;q) = \frac{1}{2}\left(1 - e^{\frac{2\rho q}{\theta}} - erf\left(\frac{-q+\rho t + v}{\sqrt{2\theta t}}\right) + e^{\frac{2\rho q}{\theta}} erf\left(\frac{q+\rho t+v}{\sqrt{2\theta t}}\right)\right). \tag{5.15}$$

The distribution function in (5.15) agrees with the transient distribution for RBM from p.49 of Harrison [9].

Returning to the general setting of Theorem 3, let

$$H_{u;\vec{x},v}(h;a) \equiv P\left(V_{u;\vec{x},v}(h) \le a\right) \tag{5.16}$$

where $V_{u;\vec{x},v}$ is the conditional net input increment defined in (5.6). If the queue length $v$ at time $u$ is large enough, then the queuing process $Q_{u;\vec{x},v}$ in (5.2) should not interact with the boundary at zero on a given interval beyond $u$; and $Q_{u;\vec{x},v} - v$ should then behave like the conditional net input process $V_{u;\vec{x},v}$ on that interval. Corollary 5 below confirms this intuition for large $v$. Because $Q_{u;\vec{x},v}(h) - v \le a \Leftrightarrow Q_{u;\vec{x},v}(h) \le v + a$, the cumulative distribution function for $Q_{u;\vec{x},v}(h) - v$ is $F_{u;\vec{x},v}(h; v+a)$.

*Corollary 5 (limiting distribution for large* $v$*): For* $F_{u;\vec{x},v}(\cdot;\cdot)$ *defined in (5.3) and* $H_{u;\vec{x},v}(\cdot;\cdot)$ *in (5.16),*

$$\lim_{v\to\infty} d_a F_{u;\vec{x},v}(h; v+a) = d_a H_{u;\vec{x},v}(h;a) = \frac{da}{\sqrt{2\pi h(\theta - \tau_u h)}} e^{\frac{-\left(a-\rho_{u;\vec{x}}h\right)^2}{2h(\theta-\tau_u h)}}$$

We recognize $\frac{1}{\sqrt{2\pi h(\theta-\tau_u h)}} e^{\frac{-\left(a-\rho_{u;\vec{x}}h\right)^2}{2h(\theta-\tau_u h)}}$ as the density of a normal distribution with mean given by

$$\lim_{v\to\infty} E\left[Q_{u;\vec{x},v}(h) - v\right] = EV_{u;\vec{x},v}(h) = \rho_{u;\vec{x}}h$$

and variance by

$$\lim_{v\to\infty} Var\left[Q_{u;\vec{x},v}(h) - v\right] = Var\, V_{u;\vec{x},v}(h) = h(\theta - \tau_u h).$$

For $Q$ defined in (5.1), $0 \le u < \Delta$ , and $0 < h < w < \Delta - u$, let

$$Q_{u;\vec{x},v}^{u+w;z}(h) \equiv Q(u+h) \ \textit{conditional on}\ Q(u) = v,\ X^{(i)}(u) = x_i \text{ for } i = 1, \dots, K \text{ and } Z(u+w) - Z(u) = z. \tag{5.17}$$

The proof of the next theorem follows that of Theorem 3, but uses Proposition 3 instead of Proposition 2.

*Theorem 4: If* $Q_{u;\vec{x},v}^{u+w;z}$ *is defined as in (5.17), then*

$$\begin{aligned} F_{u;\vec{x},v}^{u+w;z}(h;q) &\equiv P\left(Q_{u;\vec{x},v}^{u+w;z}(h) \le q\right) \\ &= \frac{1}{2}\left(1 - e^{\frac{-2q(q-z)}{w\theta}} - erf\left(\frac{-q+v+zh/w}{\sqrt{2h(\theta - \theta h/w)}}\right)\right. \\ &\qquad \left. + e^{\frac{-2q(q-z)}{w\theta}} erf\left(\frac{(q+v) - (2q-z)h/w}{\sqrt{2h(\theta - \theta h/w)}}\right)\right) \end{aligned} \tag{5.18}$$

*where* $\theta$ *is defined in (4.2).*

As a consequence of Proposition 3, the expression in (5.18) depends neither on the parameter $\rho$ nor on the values $x_i$ for $i = 1, \dots, K$.

*Corollary 6 (convergence to a point mass): For* $F_{u;\vec{x},v}^{u+w;z}(\cdot;\cdot)$ *defined in (5.18),*

$$\lim_{v\to\infty} F_{u;\vec{x},v}^{u+w;z}(h; v+a) = \frac{1}{2}\left(1 + \mathrm{erf}\left(\frac{a - zh/w}{\sqrt{2h(\theta - \theta h/w)}}\right)\right)$$

*so that*

$$\lim_{h\to w} \lim_{v\to\infty} F_{u;\vec{x},v}^{u+w;z}(h; v+a) = \begin{cases} 1, & a \ge z \\ 0, & a < z. \end{cases}$$

Corollary 6 above shows that when the increment in the net input process $Z(u+w) - Z(u) = z$ on the interval $[u, u+w]$ is given and the queue length $Q(u) = v$ that is given at the start of the interval is sufficiently large (so that $Q(\cdot)$ does not interact with the boundary on the interval), then, with probability approaching one, the increment $Q(u+w) - v$ in the queue length over the interval must equal the given increment in the net input process.

The next result is analogous to Corollary 3, but, in this case, the limit describes the conditional distribution of the queue length at the end of a measurement interval.

*Corollary 7 (conditional distribution and mean at the end of the interval): For* $F_{u;\vec{x},v}^{u+w;z}(\cdot;\cdot)$ *defined in (5.18),*

$$F_{u;\vec{x},v}^{u+w;z}(w;q) \equiv \lim_{h\to w} F_{u;\vec{x},v}^{u+w;z}(h;q) = \begin{cases} 1 - e^{-2q(q-z)/(w\theta)}, & q \geq v+z \\ 0, & q < v+z. \end{cases} \tag{5.19}$$

*and*

$$\lim_{h\to w} EQ_{u;\vec{x},v}^{u+w;z}(h) = z + v + \frac{\sqrt{2\pi\theta}\left(erf\left(\frac{2\sqrt{2}v+\sqrt{2}z}{2\sqrt{\theta}}\right)-1\right)e^{z/(2\theta)}}{4}$$

The fact that the limiting distribution function in (5.19) is equal to zero when $q < v + z$ is consistent with (5.1) since the increments of $L$ in (5.1) are non-negative by definition.

For $V_u$ defined in (5.5), $0 \leq u < \Delta$ , and $0 < h < w < \Delta - u$, let

$$V_{u;\vec{x},v}^{u+w;q}(h) \equiv V_u(h) \text{ conditional on } X^{(i)}(u) = x_i \text{ for } i = 1,\ldots,K, \quad Q(u) = v, \text{ and } Q(u+w) = q.$$

and

$$H_{u;\vec{x},v}^{u+w;q}(w;z) \equiv \lim_{h\to w} P\left(V_{u;\vec{x},v}^{u+w;q}(h) \leq z\right). \tag{5.20}$$

Our final corollary is a direct consequence of Bayes' theorem.

*Corollary 8 (density of increments of the net input process conditional on past increments and on the queue lengths at the start and end of the interval): For* $H_{u;\vec{x},v}^{u+w;q}(w;\cdot)$ *defined in (5.20),*

$$d_z H_{u;\vec{x},v}^{u+w;q}(w;z) = \frac{d_q F_{u;\vec{x},v}^{u+w;z}(w;q) \cdot d_z H_{u;\vec{x},v}(w;z)}{d_q F_{u;\vec{x},v}(w;q)}$$

*where* $F_{u;\vec{x},v}(w;q)$ *is given by (5.3),* $H_{u;\vec{x},v}(w;z)$ *by (5.16), and* $F_{u;\vec{x},v}^{u+w;z}(w;q)$ *by (5.19).*

In Fendick [18], we apply Corollary 8 to estimate conditional expectation functions for the net input and lost potential output from sample queue lengths.

## Acknowledgements

We are grateful to Ward Whitt for helpful discussions on this material.